\documentclass[11pt,reqno]{article}
\usepackage{latexsym, amsmath, amssymb, amsthm, a4, epsfig}
\usepackage{graphicx}

\newtheorem{theorem}{Theorem}[section]

\newtheorem{prop}[theorem]{Proposition}

\newcommand{\pv}{\mbox{p.v.}}

\newcommand{\ds}{\displaystyle}
\newcommand{\p}{\partial}

\newcommand{\eqnref}[1]{(\ref {#1})}

\newcommand{\Cbb}{\mathbb{C}}

\newcommand{\Rbb}{\mathbb{R}}

\newcommand{\la}{\langle}
\newcommand{\ra}{\rangle}

\newcommand{\Fcal}{\mathcal{F}}

\newcommand{\Lcal}{\mathcal{L}}

\def\BI{{\bf I}}

\def\BK{{\bf K}}

\def\BR{{\bf R}}

\def\BT{{\bf T}}

\newcommand{\Ga}{\alpha}
\newcommand{\Gb}{\beta}
\newcommand{\Gd}{\delta}

\newcommand{\Gvf}{\varphi}

\newcommand{\Gc}{\chi}

\newcommand{\Gl}{\lambda}

\newcommand{\Gm}{\mu}

\newcommand{\Gs}{\sigma}

\newcommand{\GG}{\Gamma}

\newcommand{\GO}{\Omega}

\newcommand{\BGG}{{\bf \GG}}

\newcommand{\beq}{\begin{equation}}
\newcommand{\eeq}{\end{equation}}

\def\ol{\overline}
\newcommand{\hatna}{\widehat{\nabla}}

\numberwithin{equation}{section}
\numberwithin{figure}{section}

\begin{document}

\title{Surface Riesz transforms and spectral property of elastic Neumann--Poincar\'e operators on less smooth domains in three dimensions\thanks{\footnotesize This work was
supported by NRF grants No. 2016R1A2B4011304 and 2017R1A4A1014735.}}

\author{Hyeonbae Kang\thanks{Department of Mathematics and Institute of Applied Mathematics, Inha University, Incheon
22212, S. Korea (hbkang@inha.ac.kr, k.goe.dai@gmail.com).} \and
Daisuke Kawagoe\footnotemark[2] }

\date{}
\maketitle

\begin{abstract}
It is proved in \cite{AKM18} that the Neumann--Poincar\'e operator for the Lam\'e system of linear elasticity is polynomially compact and, as a consequence, that its spectrum consists of three non-empty sequences of eigenvalues accumulating to certain numbers determined by Lam\'e parameters, if the boundary of the domain where the operator is defined is $C^\infty$-smooth. We extend this result to less smooth boundaries, namely, $C^{1,\Ga}$-smooth boundaries for some $\Ga>0$. The results are obtained by proving certain identities for surface Riesz transforms, which are singular integral operators of non-convolution type, defined by the matrix tensor on a given surface.
\end{abstract}




\section{Introduction}

The purpose of this paper is to prove certain identities for surface Riesz transforms on the boundary of a bounded domain in $\Rbb^3$, where the boundary is assumed to be $C^{1,\Ga}$ for some $\Ga>0$. We then use such identities to show that the elastic Neumann--Poincar\'e operator (the Neumann--Poincar\'e operator for the Lam\'e system of linear elasticity, abbreviated by eNP operator) on the boundary is polynomially compact. As a consequence, we show that the spectrum of the eNP operator consists of three non-empty sequences of eigenvalues accumulating to certain numbers determined by Lam\'e parameters.

Let $G(u)=(g_{ij})_{i,j=1,2}$ be a positive-definite symmetric matrix valued function on $\Rbb^2$ such that $G(u)=I$ (the identity matrix) for $u$ outside a compact set. We assume that $G$ is $C^\Ga$-smooth for some $\Ga>0$. In fact, $G$ is a metric tensor corresponding to a $C^{1,\Ga}$-smooth boundary $\p\GO$ of a certain bounded domain $\GO$ in $\Rbb^3$ (see (\ref{eq:g}) and (\ref{eq:G}) in section \ref{sec:Riesz}). Let
\beq \label{Luuv}
L(u, u-v) = \la u-v, G(u) (u-v) \ra^{-3/2}.
\eeq 
The surface Riesz transform is defined by
\beq\label{Rgdef}
R_j^g[f](u) =\frac{1}{2\pi} \text{p.v.} \int_{\Rbb^2} {(u_j-v_j)}L(u, u-v) f(v) dv, \quad j=1,2.
\eeq
Here, p.v. stands for the Cauchy principal value and $u_j$ is the $j$-th component of the point $u$. The operator $R_j^g$ is a singular integral operator of non-convolution type and bounded on $L^2(\Rbb^2)$ (or $H^{-1/2}(\Rbb^2)$) (see, for example, \cite{Stein-book-70}). The Sobolev space $H^{-1/2}(\Rbb^2)$ is of particular interest in this paper because of its relation to the spectral theory of the eNP operator.

In this paper we prove the following theorem, for presentation of which we fix notation: $A \equiv B$ for two operators $A$ and $B$ bounded on $H^s(\Rbb^2)$ ($s = 0$ or $-1/2$) means that $A - B$ is compact on $H^s(U)$ for any bounded open set $U$.

\begin{theorem}\label{thm:Riesz}
Let $R_j^g$, $j=1,2$, be surface Riesz transforms defined by the metric tensor $G$. Suppose that $G$ is $C^\Ga$-smooth for some $\Ga>0$. Then, following identities hold:
\beq\label{commute}
{R_1^g} {R_2^g}-{R_2^g} {R_1^g} \equiv 0
\eeq
and
\beq\label{squaresum}
\det(G)(g_{11} ({R_1^g})^2 + 2g_{12} R_1^g R_2^g +g_{22} ({R_2^g})^2) \equiv -I.
\eeq
\end{theorem}

It is worth mentioning that, if the surface is flat or $G$ is the identity matrix, then surface Riesz transforms are usual Riesz transforms, i.e.,
$$
R_j^g [f](u) = R_j [f](u) = \frac{1}{2\pi} \text{p.v.} \int_{\Rbb^2} \frac{u_j-v_j}{|u-v|^3} f(v) \, dv, \quad j=1,2,
$$
which are singular integral operators of convolution type, and identities (\ref{commute}) and (\ref{squaresum}) are reduced to the following ones
\beq \label{RieszId}
R_1 R_2-R_2 R_1=0 \quad\mbox{and}\quad R_1^2 + R_2^2 =-I,
\eeq
which are also proved by taking the Fourier transform. See, for example, \cite{SW-book}.

The surface Riesz transform is closely related to the eNP operator in three dimensions like the Hilbert transform is related to it in two dimensions. In fact, we show the following theorems using Theorem \ref{thm:Riesz}.

\begin{theorem}[Polynomial compactness]\label{Polynomially compact}
Let $\GO$ be a bounded domain in $\Rbb^3$ with the $C^{1,\Ga}$-smooth boundary for some $\Ga>0$. Let $\BK$ be the eNP operator on $\p\GO$ corresponding to the pair of Lam\'e parameters $(\Gl, \Gm)$. Let $p_3(t)=t(t^2-k_0^2)$ where $k_0$ is given by
\beq\label{kzero}
k_0= \frac{\mu}{2(2\mu+\Gl)}.
\eeq
Then $p_3(\BK)$ is compact. Moreover, $\BK (\BK-k_0 \BI)$, $\BK (\BK + k_0 \BI)$ and $\BK^2 - k_0^2 \BI$ are \emph{not} compact.
\end{theorem}

\begin{theorem}[Spectral structure]\label{mainthm1}
Let $\GO$ and $\BK$ be as in Theorem \ref{Polynomially compact}. The spectrum of $\BK$ consists of three non-empty sequences of eigenvalues  which converge to $0$, $k_0$ and $-k_0$, respectively.
\end{theorem}

Theorem \ref{mainthm1} is a consequence of Theorem \ref{Polynomially compact} and the spectral mapping theorem which asserts that $p_3(\Gs(\BK))=\Gs(p_3(\BK))$, where $\Gs(\BK)$ denotes the spectrum of $\BK$ (see \cite{RS-book-80}). Theorem \ref{Polynomially compact} is proved in \cite{AKM18} under the assumption that $\p\GO$ is $C^\infty$-smooth. We describe below why this assumption was needed and how it is overcome in this paper, but we first make some motivational remarks.

We will be brief here and refer to \cite{AKM18} for more informative discussion on recent development on spectral theory of the Neumann--Poincar\'e operator (abbreviated by NP operator). The NP operator, sometimes called the double layer potential, is a boundary integral operator which naturally appears when solving classical boundary value problems for the Laplace operator using layer potentials. Its study goes back to C. Neumann \cite{Neumann-87} and Poincar\'e \cite{Poincare-AM-87} as the name of the operator suggests. The NP operator, which is not a self-adjoint operator on $L^2$ in general, can be realized as a self-adjoint operator by introducing a new inner product on the Sobolev space $H^{-1/2}$ \cite{KPS-ARMA-07}. If the boundary of the domain where the NP operator is defined is $C^{1,\Ga}$-smooth for some $\Ga>0$, then the NP operator is compact. So, its spectrum consists of eigenvalues converging to $0$.

However, its counterpart for the Lam\'e system, the eNP operator, is \emph{not} compact even if the boundary is smooth \cite{DKV-Duke-88}. Therefore, it was not clear how spectrum of eNP operator looked like. But, it is proved in \cite{AJKKY} that the eNP operator in two dimensions is polynomially compact if the domain is $C^{1,\Ga}$-smooth. More precisely, if we denote the eNP operator by $\BK$, then $\BK^2 - k_0^2 I$ is compact where $k_0$ the the number given by \eqnref{kzero}, and $\Gs(\BK)$ consists of two non-empty sequences of eigenvalues converging to $k_0$ and $-k_0$, respectively. The proof of this two-dimensional result cannot be extended to three dimensions since it uses the Hilbert transform, which relates the boundary values of harmonic functions with those of their conjugate harmonic functions.

In \cite{AKM18} the three-dimensional eNP operator is expressed in terms of surface Riesz transforms $R_j^g$, and identities \eqnref{commute} and \eqnref{squaresum} are proved when $\p\GO$ is $C^\infty$-smooth. In fact, the operator $R_j^g$ is realized as a classical $\psi$do (pseudo-differential operator) and its symbol is computed (see \eqnref{symbol}). Then, calculus of $\psi$do's immediately yields those two identities. For example, ${R_1^g} {R_2^g}-{R_2^g} {R_1^g}$ is a commutator of $\psi$do's and regularizing of order $1$. Theorem \ref{Polynomially compact} is proved using \eqnref{commute} and \eqnref{squaresum} when $\p\GO$ is $C^\infty$.

If $\p\GO$ is $C^{1,\Ga}$, then the metric tensor $G$ is merely $C^\Ga$, and so calculus of $\psi$do's may not be applied. In this paper we prove Theorem \ref{thm:Riesz} by directly dealing with compositions of singular integral operators of non-convolution type. We then prove Theorems \ref{Polynomially compact} and \ref{mainthm1} following the same argument as in \cite{AKM18}

The NP operators on $C^{1,\Ga}$ boundaries and Lipschitz boundaries exhibit drastically different spectrum. It is proved recently that if $\p\GO$ has corners, then the NP operator has continuous spectrum of the connected interval symmetric with respect to $0$ whose endpoints are determined by the angle of the corner \cite{PP-arXiv} (see also \cite{BZ-arXiv, HKL-17, KLY-17, PP-JAM-14}). In this regard Theorem \ref{mainthm1} and the corresponding result in two dimensions are quite interesting. If the domain has a corner, then it is expected that the eNP operator may have a continuous spectrum. However the continuous spectrum in two dimensions may not be a connected interval, since there are two accumulation points.

This paper is organized as follows. In section \ref{sec:Comm}, we introduce an approximation of compositions of surface Riesz transforms, and we prove Theorem \ref{thm:Riesz} by using the approximation. In section \ref{sec:Riesz} we review the relations between eNP operators and surface Riesz transforms, which was obtained in \cite{AKM18}, and describe how Theorem \ref{Polynomially compact} follows from Theorem \ref{thm:Riesz}.

\section{Surface Riesz transforms and proof of Theorem \ref{thm:Riesz}}\label{sec:Comm}

In what follows, we use the notation:
\beq\label{rj}
r_j(u, v) := v_j L(u, v), \quad j=1,2,
\eeq
where $L(u, v)$ is defined by \eqnref{Luuv}. Observe that
$$
R_i^g R_j^g [f](u) = \lim_{\Gd_1, \Gd_2 \downarrow 0} \dfrac{1}{4\pi^2} \int_{|u-v| > \Gd_1} r_i(u, u - v) \int_{|v-w| > \Gd_2} r_j(v, v - w) f(w)\,dw\,dv
$$
for a.e. $u$, where the limit exists either in the point-wise sense or $L^2$-sense. Define the operator $R_{ij}$ by
$$
R_{ij} [f](u) = \lim_{\Gd_1, \Gd_2 \downarrow 0} \dfrac{1}{4\pi^2} \int_{|u-v| > \Gd_1} r_i(u, u - v) \int_{|v-w| > \Gd_2} r_j(u, v - w) f(w)\,dw\,dv
$$
for a.e. $u$. We emphasize that the difference between $R_i^g R_j^g [f](u)$ and $R_{ij} [f](u)$ lies in the $r_j$ appeared in the formulas: the first one is $r_j(v, v - w)$ while the second one is $r_j(u, v - w)$.

The following proposition is the key ingredient in proving Theorem \ref{thm:Riesz}.

\begin{prop}\label{prop:CV}
If the metric tensor $G(u)$ is $C^\Ga$ for some $\Ga>0$, then
\beq\label{RiRjRij}
R_i^g R_j^g \equiv R_{ij}.
\eeq
for $i, j = 1, 2$.
\end{prop}

Let us prove Theorem \ref{thm:Riesz} first, and then give the proof of Proposition \ref{prop:CV} after that.

\medskip
\noindent{\sl Proof of Theorem \ref{thm:Riesz}}.
Thanks to \eqnref{RiRjRij}, it suffices to prove
\beq\label{commute2}
R_{12} - R_{21} = 0
\eeq
and
\beq\label{squaresum2}
\det(G)(g_{11} R_{11} + 2g_{12} R_{12} + g_{22} R_{22}) = -I.
\eeq

Note that
\begin{align*}
R_{ij}[f](u) &= \pv \dfrac{1}{4\pi^2} \int_{\Rbb^2} r_i(u, u - v)\, \pv \int_{\Rbb^2} r_j(u, v - w)f(w)\,dw\,dv \\
&= \dfrac{1}{2\pi} \int_{\Rbb^2} \Fcal[r_i(u, \cdot)](\xi) \Fcal[r_j(u, \cdot)](\xi) \Fcal[f](\xi)\,e^{\sqrt{-1} u \cdot \xi}\,d\xi
\end{align*}
for $i, j = 1, 2$,
where $\Fcal$ denotes the Fourier transform:
$$
\Fcal[f](\xi) := \dfrac{1}{2\pi} \int_{\Rbb^2} f(x) e^{-\sqrt{-1} x \cdot \xi}\,dx.
$$
Thus, \eqnref{commute2} follows immediately.

It is proved in \cite{AKM18} that
\beq \label{symbol}
\Fcal[r_i(u, \cdot)](\xi) = \dfrac{-\sqrt{-1}}{\det(G(u))^{1/2}} \dfrac{\sum_j g^{ij}(u) \xi_j}{\sqrt{\sum_{i, j} g^{ij}(u) \xi_i \xi_j}}, \quad i = 1, 2,
\eeq
where $(g^{ij})_{i, j = 1, 2}$ is the inverse metric tensor of $G$.
Thus, we have
\begin{align*}
&\det(G)(g_{11} R_{11} + 2g_{12} R_{12} + g_{22} R_{22})[f](u)\\
&= -\det(G(u)) \int_{\Rbb^2} \dfrac{1}{\sum_{i, j} g^{ij}(u) \xi_i \xi_j} \left[ \dfrac{g_{11}(u)}{\det{G(u)}} (g^{11}(u) \xi_1 + g^{12}(u) \xi_2)^2 \right.\\
& \quad + \dfrac{2g_{12}(u)}{\det{G(u)}}(g^{11}(u) \xi_1 + g^{12}(u) \xi_2)(g^{21}(u) \xi_1 + g^{22}(u) \xi_2)\\
& \quad +\left. \dfrac{g_{22}(u)}{\det{G(u)}} (g^{21}(u) \xi_1 + g^{22}(u) \xi_2)^2 \right] \Fcal[f](\xi) e^{\sqrt{-1} u \cdot \xi}\,d\xi\\
&= -\int_{\Rbb^2} \dfrac{1}{\sum_{i, j} g^{ij}(u) \xi_i \xi_j} \big[ \left( g_{11} (g^{11})^2 + 2g_{12} g^{11} g^{21} + g_{22} (g^{21})^2 \right) (u) \, \xi_1^2 \\
& \quad + 2 \left( g_{11} g^{11} g^{12} + g^{12} (g^{11}g^{22} + g^{12} g^{21} ) + g_{22} g^{21} g^{22} \right)(u) \, \xi_1 \xi_2\\
& \quad + \left( g_{11} (g^{12})^2 + 2g_{12} g^{12} g^{22} + g_{22} (g^{22})^2 \right)(u) \, \xi_2^2 \Big] \Fcal[f](\xi) e^{\sqrt{-1} u \cdot \xi}\,d\xi.
\end{align*}
Observing that $g^{11} = g_{22} / \det(G)$, $g^{12} = g^{21} = -g_{12} / \det(G)$ and $g^{22} = g_{11} / \det(G)$, we have
\begin{align*}
&g_{11} (g^{11})^2 + 2g_{12} g^{11} g^{21} + g_{22} (g^{21})^2\\
&= \det(G) (g^{22} (g^{11})^2 - 2(g^{12})^2 g^{11} + g^{11} (g^{21})^2)\\
&= \det(G) g^{11} (g^{22} g^{11} - (g^{12})^2)\\
&= g^{11}.
\end{align*}
Similarly one can show that
$$
g_{11} g^{11} g^{12} + g^{12} (g^{11}g^{22} + g^{12} g^{21} ) + g_{22} g^{21} g^{22}=  g^{12},
$$
and
$$
g_{11} (g^{12})^2 + 2g_{12} g^{12} g^{22} + g_{22} (g^{22})^2 = g^{22}.
$$
Thus, we have
\begin{align*}
&\det(G)(g_{11} R_{11} + 2g_{12} R_{12} + g_{22} R_{22})[f](u) \\
&= -\dfrac{1}{2\pi} \int_{\Rbb^2} \dfrac{g^{11}(u) \xi_1^2 + 2g^{12}(u) \xi_1 \xi_2 + g^{22}(u) \xi_2^2}{\sum_{i, j} g^{ij}(u) \xi_i \xi_j} \Fcal[f](\xi) e^{\sqrt{-1} u \cdot \xi}\,d\xi\\
&= -\dfrac{1}{2\pi} \int_{\Rbb^2} \Fcal[f](\xi) e^{\sqrt{-1} u \cdot \xi}\,d\xi\\
&= - f(u)
\end{align*}
by the Fourier inversion formula, which completes the proof. \qed

\medskip
\noindent{\sl Proof of Proposition \ref{prop:CV}}.
In this proof we use $R_k$ for the surface Riesz transform, dropping the superscript $g$ from the notation \eqnref{Rgdef}, for ease of notation.

Note that
\begin{align*}
& R_i R_j [f](u) - R_{ij} [f](u) \\
&= \lim_{\Gd_1, \Gd_2 \downarrow 0} \dfrac{1}{4\pi^2} \int_{|u - v| > \Gd_1} \int_{|v - w| > \Gd_2} r_i(u, u  - v) \left[ r_j(v, v - w) - r_j(u, v- w) \right] f(w)\,dw\,dv.
\end{align*}
By changing the order of integrations, we see that
\beq\label{orderchange}
R_i R_j [f](u) - R_{ij} [f](u)
= \lim_{\Gd_1, \Gd_2 \downarrow 0} \dfrac{1}{4\pi^2} \int_{\Rbb^2} k_{\Gd_1,\Gd_2}(u,w)  f(w)\,dw,
\eeq
where
\beq\label{defkonetwo}
k_{\Gd_1,\Gd_2}(u,w) :=
\int_{\{|u - v| > \Gd_1\} \cap \{|v - w| > \Gd_2\}} r_i(u, u  - v) \left[ r_j(v, v - w) - r_j(u, v- w) \right] \,dv .
\eeq
We will show that $k(u, w) := \lim_{\Gd_1, \Gd_2 \downarrow 0} k_{\Gd_1,\Gd_2}(u,w)$ exists and it is weakly singular, or more precisely, if $U$ is a bounded set in $\Rbb^2$, then
\beq\label{estkonetwo}
\left| k(u, w) \right| \le \dfrac{C}{|u - w|^{2 - \Gb}}, \quad u, w \in U
\eeq
for some constant $C$, where $\Gb=3\Ga / 4$. Then, \eqnref{RiRjRij} follows from \eqnref{estkonetwo} since a weakly singular integral operator is compact on $H^s(U)$. It is worth mentioning that the order of integrations and limits can be switched in \eqnref{orderchange} since the integral in \eqnref{orderchange} is absolutely convergent as is shown in the course of proving \eqnref{estkonetwo}.

Note that
\begin{align}
&\left| L(v, v - w) - L(u, v - w) \right| \nonumber \\
& = \left| \langle v - w, G(v)(v - w) \rangle^{3/2} - \langle v - w, G(u)(v - w) \rangle^{3/2} \right|
L(v, v - w) L(u, v - w) \nonumber \\
& \leq \left| \langle v - w, G(v)(v - w) \rangle^{3/4} - \langle v - w, G(u)(v - w) \rangle^{3/4} \right| \nonumber \\
& \quad \times L(v, v - w)^{1/2} L(u, v - w)^{1/2} \left( L(v, v - w)^{1/2} + L(u, v - w)^{1/2} \right). \label{1000}
\end{align}
Here, we invoke an inequality: for all $0< p < 1$,
$$
\left| | x |^p - | y|^p \right| \leq |x - y|^p, \quad x, y \in \Rbb.
$$
So, we have
\begin{align}
& \left| \langle v - w, G(v)(v - w) \rangle^{3/4} - \langle v - w, G(u)(v - w) \rangle^{3/4} \right| \nonumber \\
& \leq \left| \langle v - w, (G(v) - G(u))(v - w) \rangle \right|^{3/4} \nonumber  \\
& \leq \left( \sum_{i, j=1}^2 |v_i - w_i| |g_{ij}(v) - g_{ij}(u)| |v_j - w_j| \right)^{3/4}  \nonumber \\
& \leq C |u - v|^{\Gb} |v - w|^{3 / 2}, \label{1001}
\end{align}
where the last inequality holds since $G(u)$ is $C^\Ga$. Note that
there exist two positive constants $C_1$ and $C_2$ such that
$$
C_1 |v-w|^{-3} \le L(v, v - w) \le C_2 |v-w|^{-3}
$$
for all $v, w \in \Rbb^2$, and similar estimates are valid for $L(u, v - w)$ as well for all $u \in \Rbb^2$.  It then follows from \eqnref{1000} and \eqnref{1001} that
$$
\left| L(v, v - w) - L(u, v - w) \right| \leq C \dfrac{|u - v|^{\Gb}}{|v - w|^3}.
$$
Since $r_j(u, v) = v_j L(u, v)$, we also have
\beq\label{rjest}
\left| r_j(v, v - w) - r_j(u, v - w) \right| \leq C \dfrac{|u - v|^{\Gb}}{|v - w|^2}.
\eeq

We assume $|u - w| > \Gd$ for some $\Gd > 0$ and take $\Gd_1$ and $\Gd_2$ such that $2 \max \{ \Gd_1, \Gd_2 \} < \Gd$. We then decompose the domain of the integral in \eqnref{defkonetwo} into two disjoint subsets:
\beq \label{decom}
\{|u - v| > \Gd_1\} \cap \{|v - w| > \Gd_2\} = A \cup B,
\eeq
where
\begin{align*}
A &:= \{|u - v| > \Gd_1\} \cap \{ |u - v| \le |v - w| \},\\
B &:= \{|v - w| > \Gd_2\} \cap \{ |v - w| < |u - v| \}.
\end{align*}
Indeed, if there existed $v \in A \cap \{|v - w| \le \Gd_2\}$, then we would have
$$
|u - w| \le |u - v| + |v - w| \le 2|v - w| \le 2 \Gd_2 < \Gd,
$$
which contradicts the assumption that $|u - w| > \Gd$. Thus, we have $A = \{|u - v| > \Gd_1\} \cap \{|v - w| > \Gd_2\} \cap \{ |u - v| \le |v - w| \}$. In the same way, we see that $B = \{|u - v| > \Gd_1\} \cap \{|v - w| > \Gd_2\} \cap \{ |v - w| < |u - v| \}$. Thus (\ref{decom}) holds.

We write
$$
k_{\Gd_1,\Gd_2}(u,w) = I_A+I_B :=
\int_{A} + \int_B r_i(u, u  - v) \left[ r_j(v, v - w) - r_j(u, v- w) \right] \,dv .
$$

We first estimate $I_A$. Observe first that
$$
A = A_1 \cup A_2 := \{ \Gd_1 < |u - v| < |u - w| /2 \} \cup \{ |u - w| / 2 < |u - v| \le |v - w| \}.
$$
According to \eqnref{rjest} we have
\begin{align*}
|I_A| \leq \int_{A_1} + \int_{A_2} \dfrac{C}{|u - v|^{2 - \Gb} |v - w|^2}\,dv =: I_{A1} + I_{A2}.
\end{align*}
If $v \in A_1$, then
$$
|v - w| \ge |u-w| - |u-v| \ge |u-w|/2,
$$
and hence
$$
I_{A1} \le \dfrac{C}{|u - w|^2} \int_{\{ \Gd_1 < |u - v| < |u - w| /2\} } \dfrac{1}{|u - v|^{2 - \Gb}}\,dv \leq \dfrac{C}{|u - w|^{2 - \Gb}}.
$$
Here and afterwards, the constant $C$ appearing in the course of estimates may differ at each occurrence, and it is independent of $\Gd$, $\Gd_1$ and $\Gd_2$. We also have
$$
I_{A2} \le \int_{\{ |u - w| / 2 < |u - v| \}} \dfrac{C}{|u - v|^{4 - \Gb}}\,dv
\leq \dfrac{C}{|u - w|^{2 - \Gb}}.
$$
Thus we have
\beq\label{IA}
\left| I_A \right| \le \dfrac{C}{|u - w|^{2 - \Gb}}.
\eeq

We now deal with $I_B$. We decompose $B$ as
$$
B=B_1 \cup B_2 := \{ \Gd_2 < |v - w| \le |u - w|/2 \} \cup \{|u - w|/2 < |v - w| < |u - v| \},
$$
and write $I_B$ as
$$
I_B  = \int_{B_1} + \int_{B_2} r_i(u, u  - v) \left[ r_j(v, v - w) - r_j(u, v- w) \right] \,dv=: I_{B1} + I_{B2} .
$$

The integral $I_{B2}$ is easy to handle. Indeed, we have
\begin{align}
|I_{B2}| &\le \int_{\{|u - w|/2 < |v - w| < |u - v|\}} \dfrac{C}{|u - v|^{2-\Gb} |v - w|^2}\,dv \nonumber\\
& \le \int_{\{|u - w|/2 < |v - w|\}} \dfrac{C}{|v - w|^{4-\Gb}}\,dv \leq \dfrac{C}{|u - w|^{2-\Gb}} \label{IB2}.
\end{align}

The rest of the proof is devoted to estimating $I_{B1}$. We first observe that
\begin{align*}
& \int_{B_1} r_i(u, u - w) \left[ r_j(w, v - w) - r_j(u, v - w) \right] \,dv \\
& = r_i(u, u - w) \int_{B_1} \left[ r_j(w, v - w) - r_j(u, v - w) \right] \,dv = 0.
\end{align*}
Thus $I_{B1}$ can be written as
\begin{align*}
I_{B1} &= \int_{B_1} \big[ r_i(u, u - v) \left[ r_j(v, v - w) - r_j(u, v - w) \right]  \\
&\quad \quad - r_i(u, u - w) \left[ r_j(w, v - w) - r_j(u, v - w) \right] \big] \,dv .
\end{align*}
We then write the integrand as
\begin{align*}
& r_i(u, u - v) \big[ r_j(v, v - w) - r_j(u, v - w) \big]
 - r_i(u, u - w) \big[ r_j(w, v - w) - r_j(u, v - w) \big] \\
&= r_i(u, u - v) \big[ r_j(v, v - w) - r_j(w, v - w) \big] \\
& \qquad + (u_i-v_i) \big[ L(u,u-v) - L(u,u-w) \big] \big[ r_j(w, v - w) - r_j(u, v - w) \big] \\
& \qquad + (v_i-w_i) L(u,u-v) \big[ r_j(w, v - w) - r_j(u, v - w) \big].
\end{align*}
Thus we have
$$
I_{B1} = J_1 + J_2 + J_3,
$$
where $J_k$ ($k=1,2,3$) corresponds to the decomposition of the integral kernel above.

We have for $J_1$ that
\begin{align*}
|J_1| & \le \int_{B_1} |r_i(u, u - v)| \big| r_j(v, v - w) - r_j(w, v - w) \big| \,dv \\
& \le \int_{B_1} \frac{C}{|u-v|^2 |v-w|^{2-\Gb}} \,dv.
\end{align*}
Note that if $v \in B_1$, then
\beq\label{u-vu-w}
|u-v| \ge |u-w|/2.
\eeq
Thus we have
$$
|J_1| \le \frac{C}{|u-w|^2} \int_{B_1} \frac{1}{|v-w|^{2-\Gb}} \,dv \le \frac{C}{|u-w|^{2-\Gb}}.
$$

To estimate $J_2$, we observe in the same way as (\ref{1000}) that
\begin{align*}
&\left| L(u, u - v) - L(u, u - w) \right|\\
&\leq \left| \langle u - v, G(u) (u - v) \rangle^{3/4} - \langle u - w, G(u) (u - w) \rangle^{3/4} \right| \\
&\quad \times L(u, u - v)^{1/2} L(u, u - w)^{1/2} \left( L(u, u - v)^{1/2} + L(u, u - w)^{1/2} \right).
\end{align*}
One can see that
\begin{align*}
& \left| \langle u - v, G(u) (u - v) \rangle^{3/4} - \langle u - w, G(u) (u - w) \rangle^{3/4} \right| \\
& \le \left| \langle u - v, G(u) (u - v) \rangle - \langle u - w, G(u) (u - w) \rangle \right|^{3/4} \\
& \le \left| \langle w - v, G(u) (u - v) \rangle + \langle u - w, G(u) (w-v) \rangle \right|^{3/4} \\
& \le C |v-w|^{3/4} \left( |u-v|^{3/4} + |u-w|^{3/4} \right).
\end{align*}
It then follows that for $v \in B_1$,
\begin{align*}
&\left| L(u, u - v) - L(u, u - w) \right| \\
& \le  C |v-w|^{3/4} \left( |u-v|^{3/4} + |u-w|^{3/4} \right) \\
& \quad \times |u-v|^{-3/2}|u-w|^{-3/2} \left( |u-v|^{-3/2} + |u-w|^{-3/2} \right) \\
& \le  C |v-w|^{3/4} \left( |u-v|^{-3/4}|u-w|^{-3/2} + |u-v|^{-3/2}|u-w|^{-3/4} \right) \\
& \quad \times \left( |u-v|^{-3/2} + |u-w|^{-3/2} \right) \\
& \le C |v-w|^{3/4} |u-w|^{-15/4},
\end{align*}
where the last inequality follows from \eqnref{u-vu-w}. Then \eqnref{rjest} and the relation $|u - v| < |v - w|$ yield that
\begin{align*}
|J_2| \le \frac{C}{|u-w|^{15/4-\Gb}} \int_{B_1} \frac{1}{|v-w|^{1/4}} dv \le \frac{C}{|u-w|^{2-\Gb}}.
\end{align*}

Similarly, one can show that
\begin{align*}
|J_3| \le \frac{C}{|u-w|^{3-\Gb}} \int_{B_1} \frac{1}{|v-w|} dv \le \frac{C}{|u-w|^{2-\Gb}}.
\end{align*}
Thus we infer that
\beq \label{IB1}
|I_{B1}| \le \frac{C}{|u-w|^{2-\Gb}}.
\eeq

From (\ref{IA}), (\ref{IB2}) and (\ref{IB1}), we have
$$
\lim_{\Gd_1, \Gd_2 \downarrow 0} |k_{\Gd_1, \Gd_2}(u, w)| \leq \dfrac{C}{|u - w|^{2 - \Gb}}
$$
for $|u - w| > \Gd$, and since $\Gd > 0$ is arbitrary, this completes the proof.
\qed

\section{Polynomial compactness of the eNP operator}\label{sec:Riesz}

It is shown in \cite{AKM18} that the eNP operator can expressed in terms of surface Riesz transforms. In this section we review it and prove Theorem \ref{Polynomially compact} using Theorem \ref{thm:Riesz}.

Let $\GO$ be a bounded domain in $\Rbb^3$ whose boundary $\p\GO$ is $C^{1,\Ga}$-smooth for some $\Ga>0$.
Let $(\Gl, \Gm)$ be the Lam\'e parameters for $\GO$ satisfying the strong convexity condition: $\Gm > 0$ and $3\Gl + 2\Gm > 0$. The isotropic elasticity tensor $\Cbb = ( C_{ijkl} )_{i, j, k, l = 1}^3$ and the corresponding Lam\'e system $\Lcal_{\Gl,\Gm}$ are defined by
$$
C_{ijkl} := \Gl \, \Gd_{ij} \Gd_{kl} + \mu \, ( \Gd_{ik} \Gd_{jl} + \Gd_{il} \Gd_{jk} )
$$
and
$$
\Lcal_{\Gl,\Gm} u := \nabla \cdot \Cbb \hatna u = \Gm \Delta u + (\Gl + \Gm) \nabla \nabla \cdot u,
$$
where $\hatna$ denotes the symmetric gradient, namely,
$$
\hatna u := \frac{1}{2} \left( \nabla u + \nabla u^T \right) \quad (T \mbox{ for transpose}). 
$$
The corresponding conormal derivative on $\p \GO$ is defined to be
$$
\p_\nu u := (\Cbb \hatna u) n = \Gl (\nabla \cdot u) n + 2\Gm (\hatna u) n \quad \mbox{on } \p \GO,
$$
where $n$ is the outward unit normal to $\p \GO$.

Let $\BGG(x) = \left( \Gamma_{ij}(x) \right)_{i, j = 1}^3$ be the Kelvin matrix of the fundamental solution to the Lam\'{e} operator $\Lcal_{\Gl, \mu}$, namely,
$$
  \Gamma_{ij}(x) =
    - \ds \frac{\Ga_1}{4 \pi} \frac{\Gd_{ij}}{|x|} - \frac{\Ga_2}{4 \pi} \ds \frac{x_i x_j}{|x|^3}, \quad x \neq 0,
$$
where
$$
  \Ga_1 = \frac{1}{2} \left( \frac{1}{\mu} + \frac{1}{2 \mu + \Gl} \right) \quad\mbox{and}\quad \Ga_2 = \frac{1}{2} \left( \frac{1}{\mu} - \frac{1}{2 \mu + \Gl} \right).
$$
The eNP operator is defined by
$$
\BK [f] (x) := \text{p.v.} \int_{\p \GO} \p_{\nu_x} {\bf \GG} (x-y) f(y) d \Gs(y) \quad \mbox{a.e. } x \in \p \GO.
$$
Here, the conormal derivative $\p_{\nu_x}\BGG (x-y)$ of the Kelvin matrix with respect to $x$-variables is defined by
$$
\p_{\nu_x}\BGG (x-y) b = \p_{\nu_x} (\BGG (x-y) b)
$$
for any constant vector $b$ (see \cite{Kup-book-65}).

Let
$$
\BK_1(x,y) = \frac{n_x (x-y)^T - (x-y) n_x^T}{2\pi |x-y|^{3}},
$$
where $n_x$ is the outward unit normal at $x$, and let
$$
\BT [f](x):= \text{p.v.} \int_{\p\GO} \BK_1(x,y) f(y) \, d \Gs(y), \quad x \in \p\GO.
$$
It is proved in \cite{AJKKY, AKM18} that
\beq\label{BKBT}
\BK \equiv k_0 \BT.
\eeq
Here \eqnref{BKBT} means that the difference $\BK - k_0 \BT$ is compact on $H^{-1/2}(\p \GO)^3$. We emphasize that $\BT$ is a singular integral operator and bounded on $H^{-1/2}(\p\GO)^3$ as well as on $L^2(\p\GO)^3$ (see \cite{CMM-AM-82}).

Denoting $n_x=(n_1(x), n_2(x), n_3(x))^T$, we have
$$
\BK_1(x,y) = \frac{1}{2\pi |x-y|^{3}}
\begin{bmatrix}
0 & K_{12} (x, y) & K_{13} (x, y) \\
- K_{12} (x, y) & 0 & K_{23} (x, y) \\
- K_{13} (x, y) & - K_{23} (x, y) & 0
\end{bmatrix},
$$
where
\begin{align*}
K_{12} (x, y) &= n_1(x)(x_2-y_2) - n_2(x) (x_1-y_1), \\
K_{13} (x, y) &= n_1(x)(x_3-y_3) - n_3(x) (x_1-y_1), \\
K_{23} (x, y) &= n_2(x)(x_3-y_3) - n_3(x) (x_2-y_2) .
\end{align*}
Let
$$
T_{ij}[f] (x) := \text{p.v.} \int_{\p\GO} \frac{K_{ij} (x, y)}{2\pi |x-y|^{3}} f(y) d\Gs(y),
$$
so that
$$
\BT =
\begin{bmatrix}
0 & T_{12} & T_{13} \\
- T_{12} & 0 & T_{23} \\
- T_{13} & - T_{23} & 0
\end{bmatrix}.
$$

Let $U$ be a coordinate chart in $\p\GO$ so that there is an open set $D$ in $\Rbb^2$ and a parametrization $\Phi: D \to U$, namely,
$$
x= \Phi(u)= (\Gvf_1(u), \Gvf_2(u), \Gvf_3(u)), \quad x \in U, \ \ u \in D.
$$
Then the metric tensor of the surface, denoted by $G(u)= (g_{ij}(u))_{i,j=1}^2$, is given by
\begin{align*}
dx_1^2+dx_2^2+dx_3^2
=g_{11} du_1^2 + 2 g_{12} du_1 du_2 + g_{22} du_2^2,
\end{align*}
where
\beq \label{eq:g}
g_{11}= |\p_1 \Phi|^2, \quad g_{12}=g_{21}=\p_1\Phi \cdot \p_2 \Phi, \quad g_{22}=|\p_2 \Phi|^2.
\eeq
Here and afterwards, $\p_j$ denotes the $j$-th partial derivative. In short, we have
\beq \label{eq:G}
G(u)= D\Phi(u)^T D\Phi(u),
\eeq
where $D\Phi$ is the $3 \times 2$ Jacobian matrix of $\Phi$. We then extend $G(u)$ to $\Rbb^2$ in such a way that $G(u)=I$ for $u$ outside a compact set. With this metric tensor, the surface Riesz transform is defined by \eqnref{Rgdef}.

Choose open sets $U_j$ ($j=1,2$) in $\p\GO$ so that $\ol{U_1} \subset U_2$ and $\ol{U_2} \subset U$. Let $\chi_j$ ($j=1,2$) be $C^{1, \Ga}$-smooth functions such that $\chi_1=1$ in $U_1$, $\mbox{supp} (\chi_1) \subset U_2$, $\chi_2=1$ in $U_2$, and $\mbox{supp} (\chi_2) \subset U$. We denote by $M_j$ the multiplication operator by $\chi_j$, i.e.,
$$
M_j [f](x) = \chi_j(x) f(x),
$$
and by $\widetilde{M}_j$ the multiplication operator by $\chi_j(\Phi(u))$ for $j=1,2$. Let $\Phi^*$ be a pull back operator, namely,
$$
\Phi^*[f](u): = f(\Phi(u)) |\p_1 \Phi \times \p_2 \Phi(u)|.
$$

For ease of notation, we set
\begin{align}
m_{11} &:= (g_{11} \p_2 \Gvf_3  - g_{12} \p_1 \Gvf_3), \label{m11} \\
m_{12} &:= (g_{21} \p_2 \Gvf_3  - g_{22} \p_1 \Gvf_3),\\
m_{21} &:= -(g_{11} \p_2 \Gvf_2  - g_{12} \p_1 \Gvf_2),\\
m_{22} &:= -(g_{21} \p_2 \Gvf_2  - g_{22} \p_1 \Gvf_2),\\
m_{31} &:= (g_{11} \p_2 \Gvf_1 - g_{12} \p_1 \Gvf_1),\\
m_{32} &:= (g_{21} \p_2 \Gvf_1 - g_{22} \p_1 \Gvf_1), \label{m32}
\end{align}
and denote by $M_{ij}$ the multiplication operator by $m_{ij}$. We emphasize that $m_{ij}$ are $C^\Ga$.

Let
\begin{align*}
X_{12} &:= \widetilde{M_2}(M_{11} {R_1^g} + M_{12}{R_2^g})\widetilde{M_1}, \\
X_{13} &:= \widetilde{M_2}(M_{21}{R_1^g}+ M_{22}{R_2^g})\widetilde{M_1} , \\
X_{23} &:= \widetilde{M_2}(M_{31}{R_1^g}+ M_{32}{R_2^g})\widetilde{M_1} ,
\end{align*}
and let
$$
\BR :=\begin{bmatrix}
	      0    &  X_{12} & X_{13}\\
	-X_{12} &      0    & X_{23}\\
	-X_{13} & -X_{23} &    0
	\end{bmatrix}.
$$
Then it is proved in \cite{AKM18} that the following relation holds:
$$
\Phi^* M_2 \BT M_1 \equiv \BR \Phi^*.
$$

Note that the crux of the matter in Theorem \ref{Polynomially compact} is that
\beq\label{p3BK}
p_3(\BK)=\BK(\BK^2-k_0^2 \BI)\equiv 0.
\eeq
In view of \eqnref{BKBT} this fact follows once we have
$$
\BT^3-\BT \equiv 0,
$$
which in turn follows from the following proposition:
\begin{prop} \label{prop:COMBR}
It holds that
\beq\label{BR3}
\BR^3 - \widetilde{M_1} \BR \equiv 0.
\eeq
\end{prop}
We refer to \cite[Section 5]{AKM18} for detailed argument to prove \eqnref{p3BK} from \eqnref{BR3}.
We now briefly show how Proposition \ref{prop:COMBR} is proved using Theorem \ref{thm:Riesz}.

\medskip
\noindent{\sl Proof of Proposition \ref{prop:COMBR}}. We first see that the following commutator relations hold:
\beq\label{3.1}
\widetilde{M_1} M_{ij} R_k^g \equiv R_k^g \widetilde{M_1} M_{ij}.
\eeq
Indeed, we have
\begin{align*}
& R_k^g \widetilde{M_1} M_{ij}[f](u) - \widetilde{M_1} M_{ij} R_k^g[f](u) \\
&= \dfrac{1}{2\pi} \int_{\Rbb^2} r_k(u, u - v) (\Gc_1(v) m_{ij}(v) - \Gc_1(u) m_{ij}(u)) f(v)\,dv,
\end{align*}
where $r_k(u, u - v)$ is defined by \eqnref{rj}. Since $m_{ij}$ is $C^{\Ga}$, we have
$$
\left| r_k(u, u - v) (\Gc_1(v) m_{ij}(v) - \Gc_1(u) m_{ij}(u)) \right| \le C|u-v|^{-2+\Ga}
$$
for some constant $C$. So \eqnref{3.1} follows, that is, $\widetilde{M_1} M_{ij} R_k^g - R_k^g \widetilde{M_1} M_{ij}$ is compact.

We then show that
\beq\label{Xij}
X_{12} X_{13} \equiv X_{13} X_{12}, \quad X_{12} X_{23} \equiv X_{23} X_{12}, \quad X_{13} X_{23} \equiv X_{23} X_{13}.
\eeq
In fact, we have
\begin{align*}
X_{12}X_{13} - X_{12}X_{13}
&= \widetilde{M_2}(M_{11} {R_1^g}+ M_{12} {R_2^g}) \widetilde{M_1}(M_{21} {R_1^g}+ M_{22} {R_2^g})\widetilde{M_1}\\
& \quad - \widetilde{M_2}(M_{21} {R_1^g}+ M_{22} {R_2^g})\widetilde{M_1}(M_{11} {R_1^g}+ M_{12} {R_2^g})\widetilde{M_1}.
\end{align*}
Here, we used the obvious identity: $\Gc_1 \Gc_2 = \Gc_1$. We then obtain using \eqnref{3.1} that
\begin{align*}
X_{12}X_{13} - X_{12}X_{13} \equiv (M_{12}M_{21}-M_{11}M_{22}) \widetilde{M_1}(R_2^g R_1^g - R_1^g R_2^g) \widetilde{M_1}.
\end{align*}
But, \eqnref{commute} implies that $\widetilde{M_1}(R_2^g R_1^g - R_1^g R_2^g) \widetilde{M_1}$ is compact. This proves the first identity in \eqnref{Xij}. The other identities there can be proved in the same way.

Cayley-Hamilton theorem and \eqnref{Xij} yield that
\beq\label{BR4}
\BR^3 + (X_{12}^2 + X_{13}^2 + X_{23}^2) \BR \equiv 0.
\eeq

One can show as before that
\begin{align*}
X_{12}^2 &= (\widetilde{M_2} M_{11} {R_1^g} \widetilde{M_1} + \widetilde{M_2} M_{12} {R_2^g} \widetilde{M_1})^2 \\
& \equiv \widetilde{M_1}(M_{11}^2 (R_1^g)^2 + 2 M_{11} M_{12} {R_1^g}{R_2^g} + M_{12}^2 (R_2^g)^2) \widetilde{M_1},
\end{align*}
likewise,
\begin{align*}
X_{13}^2 \equiv \widetilde{M_1}(M_{21}^2 (R_1^g)^2 + 2 M_{21} M_{22} {R_1^g}{R_2^g} + M_{22}^2 (R_2^g)^2) \widetilde{M_1},
\end{align*}
and
\begin{align*}
X_{23}^2 \equiv \widetilde{M_1}(M_{31}^2 (R_1^g)^2 + 2 M_{31} M_{32} {R_1^g}{R_2^g} + M_{32}^2 (R_2^g)^2)\widetilde{M_1}.
\end{align*}
Thus, we have
\begin{align*}
& X_{12}^2 + X_{13}^2 + X_{23}^2 \\
& \equiv \widetilde{M_1} \left[ \left( \sum_{i=1}^3 M_{i1}^2 \right) (R_1^g)^2 + 2 \left( \sum_{i=1}^3 M_{i1} M_{i2} \right) {R_1^g}{R_2^g} + \left( \sum_{i=1}^3 M_{i2}^2 \right) (R_2^g)^2 \right] \widetilde{M_1} .
\end{align*}
Then using the formulas \eqnref{m11}-\eqnref{m32} for $m_{ij}$ one can show that
\begin{align*}
X_{12}^2 + X_{13}^2 + X_{23}^2 &\equiv
\widetilde{M_1} \big[ (g_{22} g_{11}^2 - 2 g_{11} g_{12}^2 + g_{11} g_{12}^2) (R_1^g)^2\\
& \quad +2 (g_{22} g_{11} g_{21} - g_{11} g_{12} g_{22} - g_{21}^3 + g_{11} g_{12} g_{22}) {R_1^g}{R_2^g}\\
& \quad + (g_{22} g_{21}^2 - 2 g_{21}^2 g_{22} + g_{21} g_{22}^2) (R_2^g)^2 \big] \widetilde{M_1} \\
&= \widetilde{M_1} \det(G) (g_{11} (R_1^g)^2 + 2 g_{12} {R_1^g} {R_2^g} + g_{22} (R_2^g)^2 ) \widetilde{M_1}.
\end{align*}
It then follows from \eqnref{squaresum} that
$$
X_{12}^2 + X_{13}^2 + X_{23}^2 \equiv - \widetilde{M_1}.
$$
Now \eqnref{BR3} follows from \eqnref{BR4}. This completes the proof. \qed


\end{document}